\documentclass[12pt]{article}
\usepackage{amsmath}
\usepackage{amssymb}
\usepackage{amscd}
\usepackage{amsthm}

\theoremstyle{plain}
\newtheorem{Thm}{Theorem}[section]
\newtheorem{Conj}[Thm]{Conjecture}
\newtheorem{Prop}[Thm]{Proposition}
\newtheorem{Cor}[Thm]{Corollary}

\theoremstyle{definition}

\newtheorem{Expl}[Thm]{Example}

\theoremstyle{remark}

\numberwithin{equation}{section}

\title{On effective non-vanishing and base-point-freeness}
\author{{\it Dedicated to the memory of Professor Kunihiko Kodaira}
\\ Yujiro Kawamata}

\date{\empty}

\begin{document}

\maketitle

The celebrated Kodaira vanishing theorem implies that the cohomology
groups $H^p(X, K_X + H)$ vanish for $p > 0$ if $X$ is a smooth projective 
variety and  $H$ is an ample divisor.
It is natural to ask when $H^0(X, K_X + H)$ does not vanish.

More generally, we consider the following problem in this article.
Let $X$ be a complete normal variety, 
$B$ an effective $\mathbb{R}$-divisor on $X$,
and $D$ a Cartier divisor on $X$.
Assume that the pair $(X, B)$ is KLT (log terminal),
$D$ is nef, and that $H = D - (K_X + B)$ is nef and big (cf. \cite{KMM} for
the terminology). 
By a generalization of the Kodaira vanishing theorem 
(\cite{KMM} Theorem 1.2.5), we have $H^p(X, mD) = 0$ for any positive
integer $m$. 
The problem is to find a condition on the integer $m$ for which
the non-vanishing $H^0(X, mD) \ne 0$ holds or moreover that 
the linear system $\vert mD \vert$ is free.
By the base point free theorem (\cite{KMM} Theorem 3.1.1), 
it is known that $\vert mD \vert$ is 
free for sufficiently large integer $m$.  
Fujita's freeness conjecture implies that it should be free if 
$m \ge \dim X + 1$.
Our prediction is that $H^0(X, D) \ne 0$ always holds (Conjecture~\ref{NV}).

In \S 1, we shall derive a logarithmic version of the semipositivity 
theorem of \cite{abel}
which is used as a fundamental tool for the later sections.
In \S 2, our problem is reduced to the case where $D$ is ample.
In the rest of the paper, we consider the problem in the case $\dim X \le 4$.
In particular, we obtain a positive answer to the conjecture 
in the case $\dim X = 2$.
We also prove an existence theorem for a 4-dimensional Fano manifold.

Let us recall the terminology.
A normal variety $X$ is said to have only {\em canonical} (resp. 
{\em terminal}) singularities, if the 
following conditions are satisfied:

(1) $K_X$ is a $\mathbb{Q}$-Cartier divisor.

(2) For any birational morphism $\mu: Y \to X$ from a normal variety,
if we write $\mu^*K_X = K_Y + B^Y$ with a $\mathbb{Q}$-divisor $B^Y$ which is 
supported on the exceptional locus of $\mu$, then all the coefficients of
$B^Y$ are non-positive (resp. negative).

Let $X$ be a normal variety and $B$ an effective $\mathbb{R}$-divisor on $X$.
The pair $(X, B)$ is said to be {\em LC (log canonical)} (resp. 
{\em KLT (log terminal)}, {\em PLT (purely log terminal)}) if the 
following conditions are satisfied:

(1) $K_X + B$ is an $\mathbb{R}$-Cartier divisor.

(2) For any birational morphism $\mu: Y \to X$ from a normal variety,
if we write $\mu^*(K_X + B) = K_Y + B^Y$, then all the coefficients of
$B^Y$ are at most $1$ (resp. strictly less than $1$, at most $1$ and 
strictly less than $1$ for exceptional divisors).

\noindent
In the case where $B$ is not necessarily effective, 
the pair $(X, B)$ satisfying (1) and (2) is called 
{\em subLC (sub-log canonical)} or {\em subKLT (sub-log terminal)}.

The pair $(X, B)$ is said to be {\em properly LC} if it is LC but not KLT.
For a KLT pair $(X, B)$ and an effective $\mathbb{R}$-Cartier divisor $B'$ on
$X$, the {\em LC threshold} is defined to be 
the supremum of real numbers $t$ such that
$(X, B + tB')$ is still LC.

If the pair $(X, B)$ is LC, and if 
an irreducible component $E_j$ of $B^Y$ has the coefficient $1$, then the 
discrete valuation of the function field $\mathbb{C}(X)$ 
corresponding to the prime divisor $E_j$ is called an 
{\em LC place}, and the image $\mu(E_j) \subset X$ an {\em LC center} for the
pair $(X, B)$.
If we consider all the LC centers for all the resolutions $\mu$ 
for the fixed pair $(X, B)$ which is properly LC at a point $x \in X$,
then there exists the {\em minimal LC center} containing $x$ 
with respect to the inclusion of subvarieties of $X$ (\cite{fujita}).

Let $X$ be a reduced equi-dimensional algebraic scheme 
and $B$ an effective $\mathbb{R}$-divisor on $X$.
The pair $(X, B)$ is said to be {\em SLC (semi-log canonical)} if the 
following conditions are satisfied:

(1) $X$ satisfies the Serre condition $S_2$, 
and has only normal crossing singularities in codimension $1$.

(2) The singular locus of 
$X$ does not contain any irreducible component of $B$.

(3) $K_X + B$ is an $\mathbb{R}$-Cartier divisor.

(4) For any birational morphism $\mu: Y \to X$ from a normal variety,
if we write $\mu^*(K_X + B) = K_Y + B^Y$, then all the coefficients of
$B^Y$ are at most $1$.

We work over $\mathbb{C}$.

\section{semipositivity}

The semipositivity theorem proved in \cite{abel} 
was used for the study of algebraic fiber spaces whose fibers have  
nonnegative Kodaira dimension.
Now its logarithmic generalization will be applied for those with
negative Kodaira dimension as well.

We start with recalling the semipositivity theorem (\cite{abel} Theorem 5)
with slightly different expression: 

\begin{Thm}\label{sp}
Let $X$ and $S$ be smooth projective varieties 
and let $f: X \to S$ be a surjective morphism.
Let $n = \dim X - \dim Y$.
Assume that there exists a normal crossing divisor $\Gamma$ on $S$
such that $f$ is smooth over $S_0 = S \setminus \Gamma$.
Then the following hold:

(1) $\mathcal{F} = f_*\mathcal{O}_X(K_{X/S})$ is a locally free sheaf,
where $K_{X/S} = K_X - f^*K_S$.

(2) Let $\pi: P = \mathbb{P}(\mathcal{F}) \to S$ be the associated 
projective space bundle, and let $P_0 = \pi^{-1}(S_0)$.
Then the tautological invertible sheaf $\mathcal{O}_P(1)$ on
$P$ has a singular hermitian metric $h$ which is smooth over $P_0$ 
and such that the curvature current $\Theta$ is semipositive and that
the corresponding multiplier ideal sheaf coincides with 
$\mathcal{O}_P$.

(3) Let $X_0 = f^{-1}(S_0)$ and $f_0 = f \vert_{X_0}$.
If the local monodromies of $R^nf_{0*}\mathbb{Q}_{X_0}$ around the
branches of $\Gamma$ are unipotent,
then the Lelong number of $\Theta$ vanishes at any point of $P$.
In particular, $\mathcal{F}$ is numerically semipositive.
If $\Theta$ is strictly positive at a point on $P_0$, then
$\mathcal{O}_P(1)$ is also big.
\end{Thm}

\begin{proof}
The hermitian metric $\tilde h$ on $\mathcal{F} \vert_{S_0}$
is defined by the integration along the fiber: for $s \in S_0$
and $u,v \in \mathcal{F}_s$,
\[
\tilde h_s(u,v) = \mathrm{const.} \int_{f^{-1}(x)}u \wedge \bar v.
\]
By \cite{G}, the curvature form of $\tilde h$ is Griffiths 
semipositive.
Hence the smooth metric $h$ on $\mathcal{O}_P(1)\vert_{P_0}$ 
induced from $\tilde h$ has semipositive curvature form as well.
Moreover, it is extended to a singular hermitian metric $h$ over $P$.
The multiplier ideal sheaf is trivial because the sections of $\mathcal{F}$
are $L^2$.
In the case (3), the growth of the metric is logarithmic.
So the Lelong number vanishes, 
and the last statements follow from the regularization of positive currents
(\cite{Demailly}).
\end{proof}

The semipositivity theorem is generalized to the logarithmic case
by the covering method:

\begin{Thm}\label{sp2}
Let $X$ and $S$ be smooth projective varieties, 
let $f: X \to S$ be a surjective morphism,
and let $B$ be an effective $\mathbb{Q}$-divisor on $X$
whose support is a normal crossing divisor and whose coefficients are
strictly less than $1$.
Assume that there exists a normal crossing divisor $\Gamma$ on $S$
such that $f$ is smooth and $\mathrm{Supp}(B)$ is relative
normal crossing over $S_0 = S \setminus \Gamma$.
Let $D$ be a Cartier divisor on $X$. 
Assume that $D \sim_{\mathbb{Q}} K_{X/S} + B$.
Then the following hold:

(1) $\mathcal{F} = f_*\mathcal{O}_X(D)$ is a locally free sheaf.

(2) Let $\pi: P = \mathbb{P}(\mathcal{F}) \to S$ be the associated 
projective space bundle, and let $P_0 = \pi^{-1}(S_0)$.
Then the tautological invertible sheaf $\mathcal{O}_P(1)$ on
$P$ has a singular hermitian metric $h$ which is smooth over $P_0$ 
and such that the curvature current $\Theta$ is semipositive and that
the corresponding multiplier ideal sheaf coincides with 
$\mathcal{O}_P$.

(3) There exists a finite surjective morphism $\sigma: S' \to S$ from a
smooth projective variety $S'$ such that $\Gamma' = \sigma^{-1}(\Gamma)$ is a 
normal crossing
divisor and satisfies the following conditions:
Let $X' \to X \times_S S'$ be a birational morphism from a smooth projective
variety which is isomorphic over $S_0$ and such that 
the union of the pull-back of
the support of $B$, the pull-back of the support of $\Gamma$ and
the exceptional locus is a normal crossing divisor.
Let $f': X' \to S'$ and $\tau: X' \to X$ be the induced morphisms.
An effective $\mathbb{Q}$-divisor $B'$ on $X'$
is defined such that its coefficients are strictly 
less than $1$ and that $R = \tau^*(K_{X/S} + B) - (K_{X'/S'} + B')$ is a 
divisor.
Let $D' = \tau^*D - R$.
Then $R$ is effective, 
and the assumptions of the theorem are satisfied by $f': X' \to S'$, $B'$ 
and $D'$.
The locally free sheaf $\mathcal{F}' = f'_*\mathcal{O}_{X'}(D')$ on $S'$ 
satisfies that $\sigma^*\mathcal{F} \supset \mathcal{F}'$.
The singular hermitian metric $h$ induces a singular hermitian metric $h'$ 
on the tautological invertible sheaf 
$\mathcal{O}_{P'}(1)$ on $P' = \mathbb{P}(\mathcal{F}')$,
and the Lelong number of the curvature current $\Theta'$ 
vanishes at any point of $P'$.
In particular, $\mathcal{F}'$ is numerically semipositive.
If $\Theta'$ is strictly positive at a point on $P'_0$, then
$\mathcal{O}_{P'}(1)$ is also big.
\end{Thm}

\begin{proof}
Let $m$ be the minimal positive number such that
$mD \sim m(K_{X/S} + B)$.
We take a rational function $h$ on $X$ such that 
$\text{div}(h) = - mD + m(K_{X/S} + B)$.
Let $\pi: Y \to X$ be the normalization of $X$ 
in the field $\mathbb{C}(X)(h^{1/m})$,
and let $\mu: Y' \to Y$ be a desingularization such that
the composite morphism $g: Y' \to S$ is smooth over $S_0$.
We have 
\[
\pi_*\mathcal{O}_Y \cong 
\bigoplus_{k=0}^{m-1}\mathcal{O}_X(- kD + kK_{X/S} + 
\llcorner kB \lrcorner).
\]
The Galois group $G \cong \mathbb{Z}/m\mathbb{Z}$ acts 
on $Y$ such that the above
direct summands of $\pi_*\mathcal{O}_Y$ are eigenspaces with eigenvalues
$\text{exp}(2\pi \sqrt{-1}k/m)$.

Since $\pi$ is etale outside the support of $B$, 
$Y$ has only rational singularities, hence $\mu_*\mathcal{O}_{Y'}(K_{Y'})
= \mathcal{O}_Y(K_Y)$.
We apply Theorem~\ref{sp} to the sheaf  
$g_*\mathcal{O}_{Y'}(K_{Y'/S}) = 
f_*\pi_*\mathcal{O}_Y(K_{Y/S})$.
By duality, we have
\[
\pi_*\mathcal{O}_Y(K_Y) \cong 
\bigoplus_{k=0}^{m-1}\mathcal{O}_X(K_X + kD
- kK_{X/S} - \llcorner kB \lrcorner).
\]
By taking $k = 1$ (we may assume that $m \ge 2$), we obtain our assertions (1)
and (2) since
$f_*\mathcal{O}_X(D - \llcorner B \lrcorner)
= f_*\mathcal{O}_X(D)$.

For (3), we use the unipotent reduction theorem 
for the local monodromies of $g$
(\cite{abel}).
\end{proof}

If the base space is $1$-dimensional, we have a simpler expression:

\begin{Cor}\label{curve}
Let $X$ be a complete normal variety, 
and $B$ an effective $\mathbb{Q}$-divisor on $X$
such that the pair $(X, B)$ is KLT.
Let $f: X \to C$ be a surjective morphism to a smooth curve.
Let $D$ be a Cartier divisor on $X$ 
such that $D \sim_{\mathbb{Q}} K_{X/C} + B$.  
Then $f_*\mathcal{O}_X(D)$ is a numerically semipositive 
locally free sheaf on $C$.
\end{Cor}

\begin{proof}
Let $\mu: X' \to X$ be a log resolution for the pair $(X, B)$,
and set $\mu^*(K_X + B) = K_{X'} + B'$.  
The coefficients of $B'$ are less than $1$ and negative 
coefficients appear only for exceptional divisors of $\mu$.
We set $B' = - B'_I + B'_F$ where $B'_I$ is an effective integral divisor and
$B'_F$ is a $\mathbb{Q}$-divisor whose coefficients belong to the
interval $(0,1)$.
Since the support of $B'_I$ is exceptional for $\mu$,
we have $\mu_*\mathcal{O}_{X'}(B'_I) = \mathcal{O}_X$.
By applying the theorem to the pair $(X', B'_F)$, we deduce that
the sheaf $f_*\mathcal{O}_X(D) = 
f_*\mu_*\mathcal{O}_{X'}(\mu^*D + B'_I)$ is numerically semipositive. 
\end{proof}

\begin{Cor}\label{curve2}
Let $X, B$ and $f: X \to C$ be as in Corollary~\ref{curve}.
Let $D$ be a Cartier divisor on $X$ 
such that $H = D - (K_{X/C} + B)$ is nef and big.  
Then $f_*\mathcal{O}_X(D)$ is a numerically semipositive 
locally free sheaf on $C$.
\end{Cor}

\begin{proof}
There exists an effective $\mathbb{Q}$-divisor $B'$ such that
$(X, B + B')$ is KLT and $D \sim_{\mathbb{Q}} K_{X/C} + B + B'$.
\end{proof}

In the case of rank one sheaf, we have a more precise result which
is not used later:

\begin{Cor}\label{rank1}
In Theorem~\ref{sp2}, 
assume that $\mathcal{F} = \mathcal{O}_S(F)$ is an invertible sheaf.
Let $E$ be an effective divisor such that $E \sim D - f^*F$. 
Let $\Delta$ be the smallest $\mathbb{Q}$-divisor supported
on $\Gamma$ such that
$(X, B - E + f^*(\Gamma - \Delta))$ is subLC over the generic points
of $\Gamma$.
Then $F - \Delta$ is nef. 
\end{Cor}

\begin{proof}
Let $\mu: X' \to X$ be the log resolution of the pair 
$(X, B + E + f^*\Gamma)$.
We have $\mu^*f^*F \sim \mu^*(D - E) \sim_{\mathbb{Q}} 
\mu^*(K_{X/S} + B - E) \sim_{\mathbb{Q}} K_{X'/S} + B'$ for some
$\mathbb{Q}$-divisor $B'$.
Then our assertion is proved in \cite{adj2} Theorem 2.
\end{proof}

\section{Reduction}

We consider the following problem:

\begin{Conj}\label{NV}
Let $X$ be a complete normal variety, 
$B$ an effective $\mathbb{R}$-divisor on $X$
such that the pair $(X, B)$ is KLT, and $D$ a Cartier divisor on $X$.
Assume that $D$ is nef, and
that $H = D - (K_X + B)$ is nef and big.  
Then $H^0(X, D) \ne 0$.
\end{Conj}

This problem was considered in \cite{A} in order to construct 
ladders on log Fano varieties.
By the generalization of the Kodaira Vanishing Theorem 
(\cite{KMM} Theorem 1.2.5), we have 
$H^p(X, D) = 0$ for any positive integer $p$.
Thus the condition $H^0(X, D) \ne 0$ is equivalent to saying that 
$\chi(X, D) \ne 0$.
Our problem is a topological question, unlike the case of
the Abundance Conjectures.

The base point free theorem says that there 
exists a positive integer $m_1$ such that the linear system
$\vert mD \vert$ is free for $m \ge m_1$.
The following reduction theorem is obtained as an application of
the base point free theorem and the semipositivity theorem
with the help of the perturbation technique.

\begin{Thm}\label{red}
In Conjecture~\ref{NV}, 
one may assume that $B$ is a $\mathbb{Q}$-divisor and that $H$ is ample.
Moreover, one may assume that $D$ is also ample if one replaces $X$ suitably.
\end{Thm}

\begin{proof}
By the Kodaira lemma, 
there exists an effecive $\mathbb{R}$-divisor $E$ such that 
$B + E$ is a $\mathbb{Q}$-divisor, the pair $(X, B + E)$ is KLT, and that
$H - E$ is ample.
Therefore, we may assume that $B$ is a $\mathbb{Q}$-divisor and 
that $H$ is ample.

By the Base Point Free Theorem, there exists a proper surjective 
morphism $\phi: X \to X'$ with connected fibers to a normal projective variety 
such that $D \sim \phi^*D'$ for an ample Cartier divisor $D'$ on $X'$.
We have $H^0(X, D) \ne 0$ if and only if $H^0(X', D') \ne 0$.
We shall show that there exists an effective $\mathbb{Q}$-divisor $B'$
on $X'$ such that $(X', B')$ is KLT and $D' - (K_{X'} + B')$ is
ample. 

Since $H$ is already assumed to be ample, 
we can write $H = H_0 + 2\phi^*H'$ with $H_0$ and $H'$ being ample 
$\mathbb{Q}$-divisors.
Since $H_0$ is ample, 
there exists an effective $\mathbb{Q}$-divisor $B_0$ such that
$B + H_0 \sim_{\mathbb{Q}} B_0$ and that $(X, B_0)$ is KLT.
We set $D_0 = K_X + B_0$. Then $D_0 \sim_{\mathbb{Q}} \phi^*D'_0$ 
for $D_0' = D' - 2H'$.

We construct birational morphisms $\mu: Y \to X$ and 
$\mu': Y' \to X'$ from smooth projective varieties such that 
$\phi \circ \mu = \mu' \circ \psi$ for a morphism $\psi: Y \to Y'$.
We write $\mu^*(K_X + B_0) \sim K_Y + E_0$.
If $\mu$ and $\mu'$ are chosen suitably, then
we may assume that the conditions of \cite{adj2}~Theorem 2 are 
satisfied for $\psi$ and $E_0$.
Then there exist $\mathbb{Q}$-divisors $E'_0$ and $M$ on $Y'$ 
such that $K_Y + E_0 \sim_{\mathbb{Q}} K_{Y'} + E_0' + M$,
$\mu'_*E_0'$ is effective, $\llcorner E_0' \lrcorner \le 0$ and
$M$ is nef.
Since $H'$ is ample, there exists a $\mathbb{Q}$-divisor 
$E' \sim_{\mathbb{Q}} E'_0 + M + \mu^{\prime *}H'$ on $Y'$ such that 
$B' = \mu'_*E'$ is effective and $\llcorner E' \lrcorner \le 0$. 
Then we have $D'_0 + H' \sim_{\mathbb{Q}} K_{X'} + B'$ and $(X', B')$ is KLT.
Since $D' \sim_{\mathbb{Q}} K_{X'} + B' + H'$, we obtain our assertion.
\end{proof}

\section{Surface case}

We have a complete answer in dimension $2$.

\begin{Thm}\label{surf}
Let $X, B$ and $D$ be as in Conjecture~\ref{NV}.
Assume that the numerical Kodaira dimension $\nu(X, D)$ is at most $2$;
namely, assume that $D^3 \equiv 0$. 
Then the following hold. 

(1) $H^0(X, D) \ne 0$.

(2) The linear system $\vert mD \vert$ is free for any integer $m$ such that
$m \ge 2$.
\end{Thm}

\begin{proof}
We may assume that $\dim X = \nu(X, D) \le 2$ by Theorem~\ref{red}.
Let $\mu: X' \to X$ be the minimal resolution of singularities.
Since $\mu^*K_X - K_{X'}$ is effective, 
we can write $\mu^*(K_X + B) = K_{X'} + B'$ with $(X', B')$ being KLT.
Therefore, we may assume that $X$ is smooth.
By Theorem~\ref{red} again, 
we may also assume that $H$ is ample and that $D$ is big.

Assume first that $\dim X = 1$. 
Then the assertions follow 
immediately from the Riemann-Roch theorem.  

We assume that $\dim X = 2$ in the following.
We prove (1). By the Riemann-Roch theorem,
$\chi(X, D) = \frac 12 D(B + H) + \chi(X, \mathcal{O}_X)$.
Thus, if $\chi(X, \mathcal{O}_X) \ge 0$, then  $\chi(X, D) > 0$.
Let us assume that $\chi(X, \mathcal{O}_X) = 1 - g < 0$.  Then there exists
a surjective morphism $f: X \to C$ to a curve of genus $g$ 
whose generic fiber is isomorphism to $\mathbb{P}^1$.
By Corollary~\ref{curve2}, 
the vector bundle $f_*\mathcal{O}_X(D - f^*K_C)$ is numerically semipositive.
Since $\mathcal{O}_X(D)$ is $f$-nef, it is $f$-generated, hence
we have a surjective homomorphism $f^*f_*\mathcal{O}_X(D - f^*K_C)
\to \mathcal{O}_X(D - f^*K_C)$, and the latter sheaf is nef.
Thus $(D - f^*K_C)(B + H) \ge 0$.
Since $f^*K_C(B + H) \ge - f^*K_C \cdot K_X = 4g - 4$, 
we have $\chi(X, D) \ge g - 1 > 0$.

In order to prove (2), 
we take a general member $Y \in \vert D \vert$ as a subscheme of $X$.
We have an exact sequence $0 \to \mathcal{O}_X((m-1)D) \to \mathcal{O}_X(mD) 
\to \mathcal{O}_Y(mD) \to 0$ and $H^1(X, (m-1)D) = 0$, 
hence it is sufficient to prove the freeness 
of $\vert \mathcal{O}_Y(mD) \vert$.
Let $\mathfrak m$ be any ideal sheaf of $\mathcal{O}_Y$ of colength $1$.
We shall prove that $H^1(Y, \mathfrak m (mD)) = 0$.
By duality, it is equivalent to 
$\text{Hom}(\mathfrak m, \omega_Y(- mD)) = 0$.
Since $\text{deg }\omega_Y$ is even and $Y(B + H) > 0$, 
we have $\text{deg }\omega_Y(- mD) \le - 2$,
and we have the desired vanishing.
\end{proof}

Our bound for the freeness in Theorem~\ref{surf} is better than 
the one given by the Fujita conjecture.  But we cannot expect similar thing
in higher dimensions:

\begin{Expl} 
(1) (Oguiso) Let $X$ be a general weighted hypersurface of degree $10$ 
in a weighted projective space $\mathbb{P}(1,1,1,2,5)$.
Then $X$ is smooth, $\dim X = 3$,  
and $K_X \sim 0$.  Let $D = H = \mathcal{O}_X(1)$.
We have $H^0(X, D) \ne 0$, and $\vert 2D \vert$ is free.
But $\vert 3D \vert$ is not free, and $\vert 4D \vert$ is not very ample.

(2) Let $d$ be an odd integer such that $d \ge 3$, 
and let $X$ be a general weighted hypersurface of degree $2d$ 
in $\mathbb{P}(1, \cdots, 1,2,d)$, where
the number of $1$'s is equal to $n = \text{dim X}$.
Then $X$ is smooth. Let $D = \mathcal{O}_X(1)$.
We have $K_X \sim (d-n-2)D$, 
$D^n = 1$ and $\vert mD \vert$ is not free if $m$ is odd and $m < d$.
For example, if $n = d - 2$, then $K_X \sim 0$ and $\vert nD \vert$ is
not free.

(3) Let $d$ be an integer such that $d \not\equiv 0$ 
($\text{mod }3$) and $d \ge 4$.
Let $X$ be a general weighted hypersurface of degree $3d$ 
in $\mathbb{P}(1, \cdots, 1,3,d)$ as in (2).
Then $X$ is smooth.  Let $D = \mathcal{O}_X(1)$.
We have $K_X \sim (2d-n-3)D$, 
$D^n = 1$, and $\vert mD \vert$ is not free if
$m \not\equiv 0$ ($\text{mod }3$) and $m < d$.
For example, $\vert 2D \vert$ is not free, and
$\vert (d - 1)D \vert$ is not free if $d \equiv 2$ ($\text{mod }3$).
\end{Expl}

\section{Minimal $3$-fold}

We have so far an affirmative answer only for minimal varieties 
in the case of dimension $3$.

\begin{Prop}\label{3-fold}
Let $X$ be a $3$-dimensional projective variety with at most 
canonical singularities, and $D$ a Cartier divisor.  
Assume that $K_X$ is nef, and $D - K_X$ is nef and big.
Then $H^0(X, D) \ne 0$.
\end{Prop}

\begin{proof}
By a crepant blowings-up, 
we may assume that $X$ has only terminal singularities.
Then we have $\chi(\mathcal{O}_X) \ge - \frac 1{24}K_X c_2$ by \cite{pg},
and $3c_2 - K_X^2$ is pseudo-effective by Miyaoka \cite{Mi} (see also 
\cite{SB}).
By the Riemann-Roch theorem, we calculate
\begin{equation*}\begin{split}
h^0(X, D) &= \frac 16 D^3 - \frac 14 D^2K_X + \frac 1{12} DK_X^2
+ \frac 1{12} Dc_2 + \chi(\mathcal{O}_X) \\
&= \frac 1{12} (2D - K_X)\{\frac 16 D^2 + \frac 23 D(D - K_X) + 
\frac 16(D - K_X)^2\} \\
&\quad + \frac 1{72}(2D - K_X)(3c_2 - K_X^2) + \frac 1{24}K_Xc_2 +
\chi(\mathcal{O}_X) \\
&> 0.
\end{split}\end{equation*}
\end{proof}

\begin{Prop}\label{CY3}
Let $X$ be a complete variety of dimension $3$ 
with at most Gorenstein canonical singularities, and
$D$ a Cartier divisor.
Assume that $K_X \sim 0$ and $D$ is ample.
Let $Y \in \vert D \vert$ be a general member whose existence is 
guaranteed by Proposition~\ref{3-fold}. 
Then the pair $(X, Y)$ is LC. 
In particular, $Y$ is SLC.
\end{Prop}

\begin{proof}
Assume that $(X, Y)$ is not LC. 
Let $c$ be the LC threshold for $(X, 0)$ so that $c < 1$ and $(X, cY)$ is 
properly LC. Let $W$ be a minimal center.
By \cite{adj2} Theorem 1, for any positive rational number $\epsilon$, 
there exists an effective
$\mathbb{Q}$-divisor $B'$ on $W$ such that 
$(K_X + cY + \epsilon D)\vert_W \sim_{\mathbb{Q}} K_W + B'$ 
and $(W, B')$ is KLT.
By the perturbation technique, we may assume that $W$ is the only LC
center for $(X, cY + \epsilon D)$ 
and there exists only one LC place $E$ above $W$ 
if we replace $c$ and $\epsilon$ suitably.

Therefore, there exists a birational morphism $\mu: Y \to X$ 
from a smooth projective variety such that we can 
write $\mu^*(K_X + cY + \epsilon D) = K_Y + E + F$, 
where the support of $E + F$ is a normal crossing divisor and 
the coefficients of $F$ are strictly less than $1$.

We consider an exact sequence
\[
0 \to \mathcal{I}_W(D) \to \mathcal{O}_X(D) \to 
\mathcal{O}_W(D\vert_W) \to 0,
\]
where $\mathcal{I}_W = \mu_*\mathcal{O}_Y(- E)$ 
is the ideal sheaf for $W$.
Since $D - (K_X + cY + \epsilon D)$ is ample, 
we have $H^p(Y, \mu^*D - E + \llcorner - F \lrcorner) = 0$ and
$R^p\mu_*\mathcal{O}_Y(\mu^*D - E + \llcorner - F \lrcorner) = 0$ for $p > 0$
by the generalization of the Kodaira vanishing theorem.
Since $\mu_*\mathcal{O}_Y(\llcorner - F \lrcorner) = \mathcal{O}_X$, we 
obtain $H^1(X, \mathcal{I}_W(D)) =  0$.
Hence the homomorphism 
$H^0(X, D) \to H^0(W, D\vert_W)$ is surjective.
We have $H^0(W, D\vert_W) \ne 0$ by Theorem~\ref{surf}.
It follows that $W$ is not contained in the base locus of 
$\vert D \vert$, a contradiction.
\end{proof}

\section{Weak log Fano varieties}

The following is proved by Ambro \cite{A}.
We shall give a shorter proof of the second part 
as an application of Theorem~\ref{surf}.

\begin{Thm}\label{coindex4}
Let $X, B$ and $D$ be as in Conjecture~\ref{NV}.
Assume that there exists a positive rational number $r$ such that
$r > \dim X - 3 \ge 0$ and $- (K_X + B) \sim_{\mathbb{Q}} rD$. 
Then the following hold.

(1) $H^0(X, D) \ne 0$.

(2) Let $Y \in \vert D \vert$ be a general member.  
Then the pair $(X, B + Y)$ is PLT.
\end{Thm}

\begin{proof}
(1) is proved in \cite{A} Lemma 2.  
We recall the proof for the convenience of the reader.
We set $n = \dim X$, $d = D^n \ge 1$, $\beta = BD^{n-1} \ge 0$,
and $p(t) = \chi(X, tD)$ for $t \in \mathbb{Z}$.
Since $p(0) = 1$ and $p(-1) = p(-2) = \cdots = p(-n+3) = 0$, 
we can write 
\begin{equation*}\begin{split}
p(t) &= \frac d{n!}t^n + \frac {\beta + dr}{2(n-1)!}t^{n-1}
+ \dots \\
&= \frac d{n!}(t+1)(t+2)\cdots (t+n-3)(t^3 + at^2 + bt + 
\frac {n(n-1)(n-2)}d
\end{split}\end{equation*} 
for some numbers $a,b$.
Hence 
\[
a = \frac {n(r-n+5)}2 -3 + \frac {\beta n}{2d}.
\]
On the other hand, we have
\[
0 \le (-1)^np(-n+2) = \frac d{n(n-1)}((n-2)^2 - a(n-2) + b) - 1.
\]
Therefore
\[
h^0(X, D) = p(1) = n - 1 + (-1)^np(-n+2) + \frac {\beta + d(r-n+3)}2
> 0.
\]

(2) We may assume that $D$ is ample by the base point free theorem.
Assume that $(X, B + Y)$ is not PLT, and let $c \le 1$ be the LC threshold
so that $(X, B + cY)$ is properly LC.
Let $W$ be a minimal center.
By \cite{adj2}, for any positive rational number $\epsilon$, 
there exists an effective
$\mathbb{Q}$-divisor $B'$ on $W$ such that 
$(K_X + B + cY + \epsilon D)\vert_W \sim_{\mathbb{Q}} K_W + B'$ and 
$(W, B')$ is KLT.
Then $- (K_W + B') \sim_{\mathbb{Q}} r'D \vert_W$ for 
$r' = r - c - \epsilon$.  Since $\epsilon$ can be arbitrarily small,
we have $r' > \text{dim }W - 3$ and $r' > -1$.

If $\text{dim }W \ge 3$, then $H^0(W, D \vert_W) \ne 0$ by (1).
Otherwise, we have also $H^0(W, D \vert_W) \ne 0$ by Theorem~\ref{surf}.
By the vanishing theorem, we have $H^1(X, \mathcal{I}_W(D)) = 0$
as in the proof of Proposition~\ref{CY3}.
From an exact sequence
\[
0 \to \mathcal{I}_W(D) \to \mathcal{O}_X(D) \to \mathcal{O}_W(D\vert_W) \to 0,
\]
it follows that $W$ is not contained in the base locus of $\vert D \vert$,
a contradiction.
\end{proof}

The following result deals with the case which is just beyond the
scope of Theorem~\ref{coindex4}.

\begin{Thm}
Let $X$ be a complete variety of dimension $4$ 
with at most Gorenstein canonical singularities.
Assume that $D \sim - K_X$ is ample.
Then the following hold.

(1) $H^0(X, D) \ne 0$.

(2) Let $Y \in \vert D \vert$ be a general member.  Then $(X, Y)$ is
PLT, hence $K_Y \sim 0$ and 
$Y$ has only Gorenstein canonical singularities.
\end{Thm}

\begin{proof}
We shall prove (1) and (2) simultaneously.
Let $m$ be the smallest positive integer such that 
$H^0(X, mD) \ne 0$. We shall derive a contradiction from $m > 1$.
We take a general member $Y \in \vert mD\vert$.

Assume first that $(X, Y)$ is PLT and $m >1$. 
Then $Y$ is Gorenstein canonical.
We have an exact sequence
\[
0 \to \mathcal{O}_X(- D) \to \mathcal{O}_X((m-1)D) \to 
\mathcal{O}_Y(K_Y) \to 0.
\]
We have $\chi(X, \mathcal{O}_X(- D)) = 1$ and 
$\chi(X, \mathcal{O}_X((m-1)D)) = 0$.
On the other hand, we have
$\chi(Y, K_Y) = \frac 1{24}K_Yc_2 \ge 0$ by \cite{Mi}, a contradiction.

Next assume that $(X, Y)$ is not PLT and $m \ge 1$. 
Let $c$ be the LC threshold so that $c \le 1$ and $(X, cY)$ is 
properly LC. Let $W$ be a minimal center.
If $\text{dim }W = 3$, then we have $c \le \frac 12$.
By \cite{adj2}, for any positive rational number $\epsilon$, 
there exists an effective
$\mathbb{Q}$-divisor $B'$ on $W$ such that 
$(K_X + cY + \epsilon D)\vert_W \sim_{\mathbb{Q}} K_W + B'$ 
and $(W, B')$ is KLT.
By the perturbation technique, we may assume that $W$ is the only
center if we replace $c$ and $\epsilon$ suitably.

We consider an exact sequence
\[
0 \to \mathcal{I}_W(mD) \to \mathcal{O}_X(mD) \to 
\mathcal{O}_W(mD\vert_W) \to 0.
\]
Since $mD - (K_X + cY + \epsilon D)$ is ample, 
we have $H^1(X, \mathcal{I}_W(mD)) = 0$ by the vanishing theorem.
Hence the homomorphism 
$H^0(X, mD) \to H^0(W, mD\vert_W)$ is surjective.
If $\text{dim }W \le 2$, then we have $H^0(W, mD\vert_W) \ne 0$
by Theorem~\ref{surf}.
We shall also prove that $H^0(W, mD\vert_W) \ne 0$ 
in the case $\text{dim }W = 3$.
Then it follows that $W$ is not contained in the base locus of 
$\vert mD \vert$, a contradiction, and (1) and (2) are proved.

Assume that $\text{dim }W = 3$.
We set $r' = cm - 1 + \epsilon$ so that $K_W + B' 
\sim_{\mathbb{Q}} r'D \vert_W$.
Let $p(t) = \chi(W, tD \vert_W)$.
If we set $d = (D \vert_W)^3 > 0$ and $\delta = B'(D \vert_W)^2 \ge 0$, 
then
\[
p(t) = \frac d6 t^3 + \frac{-r'd+\delta}4 t^2 + bt + c
\]
for some numbers $b$ and $c$ by the Riemann-Roch theorem. 
By the vanishing theorem, we have
$- p(-1) \ge 0$ and $p(m-1) \ge 0$, because $r' < m - 1$.
Then 
\begin{equation*}\begin{split}
p(m) = &\frac {(m-1)(m+1)d}3
+ \frac {(m+1)(-r'd+\delta)}4 \\
&+ \frac {- p(-1) + (m+1)p(m-1)}m > 0.
\end{split}\end{equation*} 
Thus we have $H^0(W, mD\vert_W) \ne 0$.
\end{proof}

Department of Mathematical Sciences, University of Tokyo, 

Komaba, Meguro, Tokyo, 153-8914, Japan 

kawamata@ms.u-tokyo.ac.jp

\end{document}